\documentclass[twoside,11pt,a4paper]{amsart}
\usepackage{amssymb,supertabular,a4}
\usepackage{amsmath,amsthm,amscd,amssymb,graphicx}
\input xy
\xyoption{all}
\def\Hom{\hbox{\rm Hom\,}}
\def\Ext{\hbox{\rm Ext\,}}

\def\Hom{\hbox{\rm Hom\,}}
\def\dimHom{\hbox{\rm dimHom\,}}
\def\dimExt{\hbox{\rm dimExt\,}}

\def\im{\hbox{\rm im\,}}
\def\mod{\hbox{\rm mod\,}}

\def\ndp{\hbox{\rm ndp\,}}

\def\ql{\hbox{\rm ql\,}}
\def\udim{\hbox{\rm \underline{dim}\,}}
\def\ra{{\rightarrow}}

\title[GR measures processing no direct predecessors]
{The Gabriel-Roiter measures admitting no direct predecessors over
$n$-Kronecker quivers}
\author{Bo Chen}
\date{}
\address {Universit\"at zu K\"oln\\
          Mathematisches Institut\\
                   Weyertal 86-90\\
           D-50931 K\"oln\\ Germany}

\email {mcebbchen@googlemail.com}

\thanks{The author is supported by DFG-Schwerpunktprogramm
`Representation theory'.}

\begin{document}

\maketitle

\newtheorem{theo}{Theorem}[section]
\newtheorem{defi}[theo]{Definition}
\newtheorem{lemm}[theo]{Lemma}
\newtheorem{coro}[theo]{Corollary}
\newtheorem{prop}[theo]{Proposition}
\newtheorem*{main}{Theorem}
\newtheorem{Example}{Example}
\newtheorem*{Question}{Question}

\thispagestyle{empty}

\begin{abstract} Let $Q$ be an $n$-Kronecker quiver, i.e., a quiver with two vertices,
labeled by $1$ and $2$, and $n$ arrows from $2$ to $1$. We show that
there are infinitely many Gabriel-Roiter measures admitting no
direct predecessors.
\end{abstract}

\bigskip

{\footnotesize{\it Keywords.} Direct predecessor, Gabriel-Roiter
measure, $n$-Kronecker modules.}

{\footnotesize{\it Mathematics Subject Classification}(2000).
16G20,16G70}

\section{Background, preliminaries and the main result}

Let $\Lambda$ be an artin algebra and $\mod\Lambda$ the category of
finitely generated right $\Lambda$-modules. For each
$M\in\mod\Lambda$, we denote by $|M|$ the length of $M$. The symbol
$\subset$ is used to denote proper inclusion.

We first recall the original definition of Gabriel-Roiter measure
\cite{R3,R4}. Let $\mathbb{N}_1$=$\{1,2,\ldots\}$ be the set of
natural numbers and $\mathcal{P}(\mathbb{N}_1)$ be the set of all
subsets of $\mathbb{N}_1$.  A total order on
$\mathcal{P}(\mathbb{N}_1)$ can be defined as follows: if $I$,$J$
are two different subsets of $\mathbb{N}_1$, write $I<J$ if the
smallest element in $(I\backslash J)\cup (J\backslash I)$ belongs to
J. Also we write $I\ll J$ provided $I\subset J$ and for all elements
$a\in I$, $b\in J\backslash I$, we have $a<b$. We say that $J$ {\bf
starts with} $I$ if $I=J$ or $I\ll J$. Thus $I<J<I'$ with $I'$
starts with $I$ implies that $J$ starts with $I$.

For each $M\in\mod\Lambda$, let $\mu(M)$ be the maximum of the sets
$\{|M_1|,|M_2|,\ldots, |M_t|\}$, where $M_1\subset M_2\subset \ldots
\subset M_t$ is a chain of indecomposable submodules of $M$. We call
$\mu(M)$ the {\bf Gabriel-Roiter (GR for short) measure}  of $M$. A
subset $I$ of $\mathcal{P}(\mathbb{N}_1)$ is called a GR measure if
there is an indecomposable module $M$ with $\mu(M)=I$. If $M$ is an
indecomposable $\Lambda$-module, we call an inclusion $X\subset M$
with $X$ indecomposable a {\bf GR inclusion} provided
$\mu(M)=\mu(X)\cup\{|M|\}$, thus if and only if every proper
submodule of $M$ has Gabriel-Roiter measure at most $\mu(X)$. In
this case, we call $X$ a {\bf GR submodule} of $M$.  Note that the
factor of a GR inclusion is indecomposable.

Using Gabriel-Roiter  measure, Ringel obtained a partition of the
module category for any artin algebra of infinite representation
type \cite{R3,R4}: there are infinitely many GR measures $I_i$ and
$I^i$ with $i$ natural numbers, such that
$$I_1<I_2<I_3<\ldots\quad \ldots <I^3<I^2<I^1$$ and such that any
other GR measure $I$ satisfies $I_i<I<I^j$ for all $i,j$. The GR
measures $I_i$ (resp. $I^i$) are called take-off (resp. landing)
measures. Any other GR measure is called a central measure. An
indecomposable module $M$ is  called a  take-off (resp. central,
landing) module if its GR measure $\mu(M)$ is a take-off (resp.
central, landing) measure.

Let $I,I'$ be two GR measures for $\Lambda$. We call $I'$  a {\bf
direct successor} of $I$ if, first, $I<I'$ and second, there does
not exist a GR measure $I''$ with $I<I''<I'$. The so-called {\bf
Successor Lemma} in \cite{R4} states that any GR measure $I$
different from $I^1$, the maximal one, has a direct successor.
However, there is no `Predecessor Lemma'. For example, the GR
measure of a homogeneous simple module over a tame quiver of type
$\widetilde{\mathbb{A}}_n$ has no direct predecessor \cite{Ch4}.
More general, for a tame quiver (of type $\widetilde{\mathbb{A}}_n$,
$\widetilde{\mathbb{D}}_n$, $\widetilde{\mathbb{E}}_6$,
$\widetilde{\mathbb{E}}_7$ or $\widetilde{\mathbb{E}}_8$), the
minimal central GR measure exists, which obviously has no direct
predecessor \cite{Ch3}.  Let us denote by $\ndp(\Lambda)$ the number
of the GR measures admitting no direct predecessors. It is clear
that any GR measure over a representation-finite artin algebra has a
direct predecessor, i.e., $\ndp(\Lambda)=0$. It was shown in
\cite{Ch4} that for each tame hereditary algebra over an
algebraically closed field $\ndp(\Lambda)<\infty$ , and thus $1\leq
\ndp(\Lambda)<\infty$ by above mentioned facts. Therefore,  it is
natural to ask if the number of the GR measures admitting no direct
predecessors relates to representation type of hereditary algebras.
More general, we want to know whether a finite dimensional algebra
$\Lambda$ over an algebraically closed field is of wild type implies
that there are infinitely many GR measures admitting no direct
predecessors, i.e., $\ndp(\Lambda)=\infty$, and vice versa.

Let  $Q'$ be the following quiver:
$$\xymatrix{1\ar@{->}[r] & 2\ar@/^/[r]\ar@/_/[r]& 3\\},$$ and $\mod(Q)$ the
category of finite dimensional representations over an algebraically
closed field.   In \cite{Ch5}, it was proved that the GR measure
$\mu^m=\{1,2,4,\ldots, 2m,2m+1\}$ has no direct predecessor for each
$m\geq 1$. The embedding of the $2$-Kronecker quiver into $Q'$ was
used in the discussion.

In this note, a second example of wild quiver with infinitely many
GR measures admitting no direct predecessors will be presented.
Namely the so called $n$-Kronecker quiver with $n\geq 3$ will be
studied. Using the same method in \cite{Ch5}, we show the following
theorem:

\begin{main}\label{ndp} Let $Q$ be an $n$-Kronecker quiver with $n\geq 3$.
Then for each  $m\geq 1$, $\mu^m=\{1,2,4,\ldots,2m,2m+1\}$ is a GR
measure which does not admit a  direct predecessor. Thus
$\ndp(Q)=\infty$.
\end{main}

\section{The GR measures for $n$-Kronecker quiver}

Let $Q$ be the following  $n$-Kronecker quiver:
$$\xymatrix{ 2\ar@/^1pc/[rr]^{\alpha_1}\ar@/_1pc/[rr]_{\alpha_n}&\vdots& 1\\}.$$
Note that $Q$ is of finite representation type if $n=1$ and of  tame
type if $n=2$. If $n\geq 3$, then  $Q$ is of wide representation
type.  We recall some facts of representations of quivers. The best
references are \cite{ARS,R2}. We also refer to \cite{K,R1} for
general structures of representations of wild quivers. Let $k$ be an
algebraically closed field. A representation for $Q$ over $k$ is
simply called an $n$-Kronecker module.

The
Cartan matrix and the Coxeter matrix are the following: $$C=\left(\begin{array}{cc}1 & 0 \\
n & 1\\\end{array}\right),\quad \Phi=-C^{-t}C=
\left(\begin{array}{cr}n^2-1 & n \\
                        -n & -1 \\ \end{array}\right),\quad \Phi^{-1}=
\left(\begin{array}{rc}-1 & -n \\
                         n & n^2-1 \\
                         \end{array}\right)$$
The dimension vectors can be calculated using $\udim\tau M=(\udim
M)\Phi$ if $M$ is not projective and $\udim \tau^{-1}N=(\udim
N)\Phi^{-1}$ if $N$ is not injective, where $\tau$ denotes the
Auslander-Reiten translation. The Euler form is $\langle
\underline{x},\underline{y} \rangle=x_1y_1+x_2y_2-nx_1y_2$. Then for
two indecomposable modules $X$ and $Y$,
$$\dimHom(X,Y)-\dimExt^1(X,Y)=\langle\udim X,\udim Y\rangle.$$

Now let us assume that $n\geq 2$. The Auslander-Reiten quiver of $Q$
consists of one preprojective component, one preinjective component
and infinitely many regular ones. An indecomposable regular module
$X$ is called quasi-simple if the Auslander-Reiten sequence starting
with $X$ has an indecomposable middle term. For each indecomposable
regular module $M$, there is a unique quasi-simple module $X$ and a
unique natural number $r\geq 1$ (called quasi-length of $M$ and
denoted by $\ql(M)=r$) such that there is a sequence of irreducible
monomorphisms $X=X[1]\ra X[2]\ra\ldots\ra X[r]=M$.

The preprojective component is the following (note that there are
namely $n$ arrows from $P_i$ to $P_{i+1}$):
$$\xymatrix@R=8pt@C=12pt{
   &P_2=(1,n)\ar@{->}[rd]\ar@{.}[rr] &&
    P_4=(n^2-1,n^3-2n)\ar@{->}[rd]&\ldots&
     \\
  P_1=(0,1)\ar@{->}[ru]\ar@{.}[rr] &&
  P_3=(n,n^2-1)\ar@{->}[ru]\ar@{.}[rr]&&
  P_5&\ldots\\}$$
Similarly, the preinjective component is of the following shape:
$$\xymatrix@R=8pt@C=12pt{
   &\ldots&Q_3=(n^3-2n,n^2-1)\ar@{->}[rd]\ar@{.}[rr] &&
    Q_1=(n,1)\ar@{->}[rd]&
     \\
  \ldots&Q_4\ar@{->}[ru]\ar@{.}[rr] &&
  Q_2=(n^2-1,n)\ar@{->}[ru]\ar@{.}[rr]&&
  Q_0=(1,0).\\}$$

\subsection{The partition using GR measrue}
Before describing the take-off part and  the landing part of
$n$-Kronecker quiver for $n\geq 2$.  we need some properties of GR
measures. The following property was proved in \cite{R3}:
\begin{prop}\label{GA} Let $\Lambda$ be an artin algebra and $X$ and
$Y_1,Y_2,\ldots,Y_r$ be indecomposable modules. Assume that
$X\stackrel{f}{\ra}\oplus_{i=1}^r {Y_i}$ is a monomorphism.
\begin{itemize}
  \item[(1)] $\mu(X)\leq \textrm{max}\{\mu(Y_i)\}$.
  \item[(2)] If $\textrm{max}\{\mu(Y_i)\}=\mu(X)$,
          then $f$ splits.
\end{itemize}
\end{prop}
We collect some properties of GR inclusions in the following lemma.
The proof can be found for example in \cite{Ch1,Ch2}.
\begin{lemm}\label{GR} Let $\Lambda$ be an artin algebra and $X\subset M$ a GR
inclusion.
\begin{itemize}
  \item[(1)] If all irreducible maps to $M$ are monomorphisms, then
              the GR inclusion is an irreducible map.
  \item[(2)] Every irreducible map to $M/X$ is an epimorphism.
  \item[(3)] There is an irreducible monomorphism $X\ra Y$ with $Y$
             indecomposable and an epimorphism $Y\ra M$.
\end{itemize}
\end{lemm}

Let us denote by $I_i$ (resp. $I^i$) the take-off (resp. landing)
measure obtained using GR measure and by $\mathcal{A}(I)$ the set of
isomorphism classes of indecomposable modules (or representatives of
these isomorphism classes) with GR measure $I$. Obviously,
$I_1=\{1\}$ and $\mathcal{A}(I_1)$ contains precisely all simple
modules. It is easily seen that $I_2=\{1,r\}$, where $r$ is the
largest possible length of a local module of Loewy length $2$.

\begin{prop} Let $n\geq 2$ and $Q$ be the $n$-Kronecker quiver.
\begin{itemize}
   \item[(1)]  $\mathcal{A}(I_r)=\{P_r\}$ for all $r\geq 2$.
              Thus the take-off part contains precisely
                the simple injective module and the indecomposable preprojective modules.
  \item[(2)] $\mathcal{A}(I^r)=\{Q_{r}\}$ for all $r\geq 1$.
              Thus the landing part contains precisely the non-simple indecomposable preinjective module.
   \item[(3)]  An indecomposable module is in central part if and only if
               it is a regular module.
\end{itemize}
\end{prop}

\begin{proof}
(1) By Lemma \ref{GR}(1), $P_i$ is the unique, up to isomorphism, GR
submodule of $P_{i+1}$.  We proceed by induction. For $r=2$, the
assertion holds by the description of $I_2$. Assume that
$\mu(M)=I_{r+1}$ for some indecomposable module $M$. Since $M$ is
not simple, we can assume that $Y$ is a GR submodule of $M$. Then
$\mu(Y)=I_i\leq I_r$ for some $i\leq r$, and thus $Y\cong P_i$ by
induction. It follows from Lemma \ref{GR}(3) that there is an
epimorphism $P_{i+1}\ra M$. In particular $|M|\leq |P_{i+1}|$. If
the equality does not hold, then
$$I_{r+1}=\mu(M)=I_i\cup\{|M|\}>I_i\cup\{|P_{i+1}|,\ldots,|P_r|,|P_{r+1}|\}
>I_i\cup\{|P_{i+1}|,\ldots,|P_r|\}=I_r.$$
This is a contradiction because the GR measure
$\mu(P_{r+1})=I_i\cup\{|P_{i+1}|,\ldots,|P_r|,|P_{n+1}|\}$ lies
between $I_r$ and $I_{r+1}$, and $I_{r+1}$ is a direct successor of
$I_r$. Therefore, $|P_{i+1}|=|M|$ and thus $P_{i+1}\cong M$. Since
$\mu(M)=I_{r+1}$, we have $i=n$ and thus
$\mathcal{A}(I_{r+1})=\{P_{r+1}\}$.

(2) It was proved in \cite{R3} that all modules lying in the landing
part are preinjective (in more general sense). Since there is a
short exact sequence $0\ra Q_{r+1}\ra Q_r^n\ra Q_{r-1}\ra 0 $ for
each $r\geq 1$, $\mu(Q_{r+1})<\mu(Q_r)$ by Proposition \ref{GA}.
Since landing modules are preinjective, $\mathcal{A}(I^r)$ contains
precisely one isomorphism class, namely $Q_r$.

Statement (3) is a direct consequence of (1) and (2).
\end{proof}
In general, the take-off part of a bimodule algebra can be similarly
described. We refer to \cite{R5} for details.

\subsection{The GR measures of $2$-Kronecker modules}\label{example}
As an example, we describe the GR measures for $2$-Kronecker modules
to be used in our later discussion. These follow easily from direct
calculation.     For more properties of general tame quivers of type
$\widetilde{\mathbb{A}}_n$, we refer to \cite{Ch2, Ch4}.

The GR measure of the indecomposable preprojective module with
dimension vector $(m,m+1)$ is $\{1,3,5,\ldots,2m+1\}$.

Every indecomposable regular module with dimension vector $(m,m)$
has GR measure $\{1,2,4,6,\ldots, 2m\}$.

The GR measure of the indecomposable preinjective module with
dimension vector $(m+1,m)$ is $\{1,2,4,\ldots, 2m,2m+1\}$.
\smallskip
\\
{\bf Remark.} Comparing to the main theorem, the GR measure
$\mu^m=\{1,2,4,\ldots,2m,2m+1\}$  in $2$-Kronecker case has a direct
predecessor $\mu^{m+1}$, for each $m\geq 1$.

\subsection{GR measures admitting no direct predecessors}

Now we assume that $Q$ is of wild type, i.e., $n\geq 3$. The
indecomposable regular modules $X$ with dimension vector $(1,1)$ or
$(2,1)$ are of special interest because they have no proper regular
factor modules. Thus any non-zero homomorphism from $X$ to an
indecomposable regular $M$ is a monomorphism.

\begin{lemm} \begin{itemize}
      \item[(1)] An indecomposable module $M$ has GR measure $\{1,2\}$ if and
                only if $\udim M=(1,1)$.
      \item[(2)] An indecomposable module $M$ has
               GR measure $\{1,2,3\}$ if and only if $\udim M=(2,1)$.
          \end{itemize}
\end{lemm}
\begin{proof} Straightforward.
\end{proof}

\begin{lemm}\label{a0} Let $X$ be an indecomposable module with dimension vector $(1,1)$.
Then for each $i\geq 1$ the GR measure $\mu(\tau^i X)$ starts with
$\{1,2,3\}$, i.e., $\mu(\tau^iX)=\{1,2,3,\ldots\}$.
\end{lemm}

\begin{proof}If $X$ is an indecomposable module with dimension vector
$(1,1)$.
Then $$\udim\tau X=(1,1)\left(\begin{array}{cr}n^2-1 & n \\
                        -n & -1 \\
                        \end{array}\right)=(n^2-n-1,n-1).$$
Thus $\dimHom(X,\tau X)-\dimExt^1(X,\tau X)=\langle\udim X,\udim\tau
X\rangle=n^2-n-1+n-1-n^2+n=n-2\geq 1$ for each $n\geq 3$. In
particular, $\Hom(X,\tau X)\neq 0$ and there exists a monomorphism
from $X$ to $\tau X$. Therefore, there exists a  monomorphism
$\tau^iX\ra\tau^{i+1}X$ for each $i\geq 0$. In particular,
$\mu(\tau^iX)>\mu(\tau^jX)$ for $i>j\geq 1$.

Let $Y$ be an indecomposable module with dimension vector $(2,1)$.
Then $$\dimHom(Y,\tau X)-\dimExt^1(Y,\tau X)=\langle\udim
Y,\udim\tau X\rangle=n-3\geq 0.$$ Thus if $n\geq 4$, then there is a
monomorphism $Y\ra \tau X$. In particular, $\mu(\tau^iX)\geq\mu(\tau
X)>\mu(Y)=\{1,2,3\}$ for each $i\geq 1$.

Assume that $n=3$.  We show that $\tau X$ contains an indecomposable
submodule with dimension vector $(2,1)$. Note that $\udim \tau
X=(5,2)$.  Since each indecomposable module with dimension vector
$(2,1)$ has no proper regular factor, it is sufficient to show that
there is an indecomposable module $Y$ with dimension vector $(2,1)$
such that $\Hom(Y,\tau X)\cong \mathbb{D}\Ext^1(X,Y)\cong
\Hom(\tau^{-1}Y,X)\neq 0$. Thus it is sufficient to show that there
is an indecomposable module $Y'$ with dimension vector
$(1,2)=(2,1)\Phi^{-1}$ such that $\Hom(Y',X)\neq 0$. Let $P_1$ be
the simple projective module. Then $\Ext^1(X,S_1)\cong
\mathbb{D}\Hom(S_1,\tau X)\neq 0$. Thus there is a non-split short
exact sequence $0\ra S_1\ra Y\ra X\ra 0.$ It is easily seen that the
middle term $Y'$ is indecomposable and has dimension vector $(1,2)$.
Since $\tau X$ contains an indecomposable submodule with dimension
vector $(2,1)$, we have  $\mu(\tau ^iX)\geq\mu(\tau X)>\{1,2,3\}$
for all $i\geq 1$.
\end{proof}

\begin{lemm}\label{a} For each $m\geq 1$, let $\mu_m=\{1,2,4,\ldots, 2m\}$.
\begin{itemize}
   \item[(1)] $\mu_m$ is a GR measure and if $M$ is an indecomposable module
             with $\mu(M)=\mu_m$, then $\udim M=(m,m)$.
    \item[(2)] If $\mu(M)=\mu_m$, then each indecomposable
             regular factor module of $M$ contains some indecomposable submodule
              with dimension vector $(1,1)$.
\end{itemize}
\end{lemm}

\begin{proof}
(1) Let $M$ be an indecomposable  $2$-Kronecker module with
dimension vector $(m,m)$. Then $M$ is  an  indecomposable
$n$-Kronecker module as well and  $M$ has GR measure $\mu_m$
(Section \ref{example}). Thus $\mu_m$ is a GR measure. Now assume
that $M$ is an indecomposable module with GR measure $\mu(M)=\mu_m$.
If $m=1$, then $\mu(M)=\{1,2\}$ and thus $\udim M=(1,1)$. Assume
that $m>1$. A GR submodule $N$ of $M$ has GR measure $\mu_{m-1}$. By
induction, $\udim N=(m-1,m-1)$. Since a GR factor of length $2$ has
dimension vector $(1,1)$, $\udim M=(m,m)$.

(2) We use induction on $m$. This is clear for $m=1$. If $m=2$, then
the dimension vector of a proper indecomposable regular factor is
$(1,1)$ or $(2,1)$. Now assume that $m\geq 2$. Let $Y$ be a proper
regular factor of $M$ and $K$ be the kernel of the canonical
projective. Let $N$ be a GR submodule of $M$. Then we have the
following diagram:
$$\xymatrix{&&N\ar[d]^f&&\\0\ar[r] & K\ar[r]^{\iota} & M\ar[r]^{\pi} & Y\ar[r] & 0\\ }$$
If the composition $\pi f$ is not zero, then the image of $\pi f$ is
a regular factor of $N$.  Therefore  by induction, $\im \pi f$, and
thus $Y$, contains a submodule with dimension vector $(1,1)$. If the
composition $\pi f$ is zero, then $l$ factors through $\iota$ by
definition of kernel.  In particular, there is a monomorphism
$N\stackrel{g}{\ra}K$. Since $N$ is a GR submodule of $M$, it
follows that $N$ is isomorphic to a direct summand of $K$
(Proposition \ref{GA}). Since $Y$ is not simple and $|K|\geq |X|$,
we have $K\cong N$. Thus $|Y|=2$ and $\udim Y=(1,1)$.
\end{proof}

\begin{lemm}\label{a1} Let $\mu^m=\{1,2,4,\ldots,2m,2m+1\}$ for each $m\geq
1$. Then
  \begin{itemize}
     \item[(1)] $\mu^m$ is a GR measure.
     \item[(2)] If $M$ is an indecomposable module with GR measure $\mu(M)=\mu^m$, then $\udim M=
                (m+1,m)$.
   \end{itemize}
\end{lemm}

\begin{proof} It is known that each indecomposable $2$-Kronecker module with dimension vector $(m+1,m)$
has GR measure $\mu^m$. Thus $\mu^m$ is a GR measure for
$n$-Kroneker quiver. We have seen that an indecomposable module $M$
has GR measure $\{1,2,4,\ldots,2m\}$ implies that $\udim M=(m,m)$.
Thus an indecomposable module with GR measure $\mu^m$ has dimension
vector $(m+1,m)$.
\end{proof}

\begin{coro}\label{a2} Let $M$ be an indecomposable module with GR measure $\mu(M)=\mu^m$. Then
        each indecomposable regular factor module of $M$ contains some indecomposable submodule with
    dimension vector $(1,1)$.
\end{coro}

\begin{proof} By Lemma \ref{a1}, $\udim M= (m+1,m)$. we have a
short exact sequence
$$0\ra N\stackrel{\iota}{\ra} M\ra M/N\ra 0$$
where  $\iota$ is a GR inclusion. Thus $\udim N=(m,m)$ and  the
factor $M/N$ is an injective simple module.  Let
$M\stackrel{\pi}{\ra}Y$ be an epimorphism with $Y$ an indecomposable
regular module.  Then $\Hom(M/N,Y)=0$. By definition of cokernel,
the composition $N\stackrel{\pi\iota}{\ra}Y$ is not zero. Therefore,
by Lemma \ref{a}, the image of $\pi \iota$, and thus $Y$, contains a
submodule with dimension vector $(1,1)$.
\end{proof}

\begin{lemm}\label{b} Let $m\geq 1$ and $M$ be an indecomposable module such that
$\mu(M)>\mu^m$.  Then $\mu(M)$ starts with
$\mu(M)=\{1,2,\ldots,2t,2t+1\}$ for some $1\leq t\leq m$. In
particular, $M$ contains an indecomposable submodule with GR measure
$\mu^t$.
\end{lemm}

\begin{proof}
This follows directly from the definition of GR measure.
\end{proof}

\begin{lemm}\label{c} If $M$ is an indecomposable module such that
$\mu=\mu(M)$ is a direct predecessor of $\mu^m$ for some $m$. Then
$M$ is regular and $|M|>2m+1$.
\end{lemm}

\begin{proof}
Note that a direct predecessor of a central measure is a central
measure as well. Thus $M$ is regular since $\mu^m$, and thus
$\mu(M)$, is a central measure. Since $\mu(M)$ is a direct
predecessor of $\mu^m$, we have $\mu_m\leq \mu(M)$. However, the
equality does not hold because $\mu_m<\mu_{m+1}<\mu^m$.  Thus
$\mu_m<\mu(M)<\mu^m$ and $|M|>2m+1$.
\end{proof}

{\it proof of Theorem.} For the purpose of a contradiction, we
assume that $M$ is an indecomposable module such that $\mu(M)$ is a
direct predecessor of $\mu^m$ for a fixed $m\geq 1$. Thus $M$ is a
regular module by Lemma \ref{c} and $M\cong X[r]$ for some
quasi-simple $X$ and $r\geq 1$. Again by Lemma \ref{c}, we have
$|X[r+1]|>|X[r]|>2m+1$.   It follows that
$\mu(M)=\mu(X[r])<\mu^m<\mu(X[r+1])$. Thus $X[r+1]$ contains a
submodule $Y$ with GR measure $\mu^t$ for some $1\leq t\leq m$
(Lemma \ref{b}). Note that $\udim Y=(t+1,t)$ and $\mu(Y)\geq \mu^m$.
We claim that $\Hom(Y,\tau^{-r}X)=0$. If this is not the case, then
by Corollary \ref{a2}, the image of a nonzero homomorphism, in
particular $\tau^{-r}X$, contains a submodule $Z$ with dimension
vector $(1,1)$. Therefore, there is a monomorphism $\tau^r Z\ra X$,
and thus, by Lemma \ref{a0},
 $$\mu^m>\mu(M)=\mu(X[r])\geq \mu(X)>\mu(\tau^r Z)> \{1,2,3\},$$ which is a contradiction.
Since there is a short exact sequence $$0\ra X[r]\ra X[r+1]\ra
\tau^{-r}X\ra 0$$ and $\Hom(Y,\tau^{-r}X)=0$, the inclusion $Y\ra
X[r+1]$ factors through $X[r]$. Therefore, there is a monomorphism
$Y\ra X[r]$. It follows that
$$\mu(X[r])\geq \mu(Y)=\mu^t \geq \mu^m>\mu(M).$$  This contradiction implies that $\mu^m$ has no direct
predecessor for each $m\geq 1$. In particular, $\ndp(Q)=\infty$.
\qed

\end{document}